\documentclass[10pt,onecolumn]{article}
\usepackage{amsfonts,amsmath,amsthm,epsfig,verbatim,amssymb,amscd,eucal,bbm}
\def\P{{\mbox{\bf P}}}

\def\PP{{\mathbb{P}}}
\def\EE{{\mathbb{E}}}

\def\Ints{{\hbox{$\it Z\hskip-3.6pt Z$}}}
\def\Nats{{\mathbb N}}
\def\Ints{{\mathbb Z}}

\def\nv{n} %number of vertices - decide N or n
\def\ne{E} %2\ne is the number of edges - E or M or m
\def\dvec{{\bf d}}
\def\dmean{{\overline d}}
\def\whp{{\bf whp}}

\newtheorem{theorem}{Theorem}
\newtheorem{lemma}{Lemma}
\newtheorem{prop}{Proposition}
\begin{document}

\title{Exponential random graphs as models of overlay networks}

\author{ M. {\sc Draief} \thanks{Imperial College
E-mail: {\tt m.draief@imperial.ac.uk}},~A. {\sc Ganesh}
\thanks{University of Bristol, E-mail: {\tt A.Ganesh@bristol.ac.uk} } ~ and ~ L. {\sc Massouli\'e}
\thanks{Thomson Research E-mail: {\tt laurent.massoulie@thomson.net} }}

\date{}

\maketitle

\begin{abstract}
In this paper, we give an analytic solution for graphs with $n$ nodes and $E$ edges for which the probability of obtaining a given graph $G$ is $\mu(G)=e^{-\beta\sum_{i=1}d_i^2}$, wherer $d_i$ is the degree of node $i$. We describe how this model naturally appears in the context of load balancing in communication networks, namely Peer-to-Peer overlays. We then analyse the degree distribution of such graphs and show that the degrees are concentrated around their mean value. Finally, we derive asymptotic results on the number of edges crossing a graph cut and use these results  $(i)$ to compute the graph expansion and conductance, and $(ii)$ to analyse the graph resilience to random failures.
\end{abstract}



\section{Introduction}

Random graphs provide a way of modelling large and complex networks, and of studying stochastic processes on such networks. Early work on this topic goes back to the famous random graph or Bernoulli graph introduced by Solomonoff and Rapoport \cite{SoRa51} in the early 1950s and studied by Erd\"os-R\'enyi \cite{ER60} a decade later. The Bernoulli random graph model is, however, rather simplistic and fails to capture important features of many real-world networks. This has stimulated work on a number of other random graph models.
Exponential random graphs were first introduced in the early 1980s by Holland and Leinhardt \cite{HL81} based on the work of Besag \cite{Besag74}. More recently Frank and Strauss \cite{FS86} studied a subclass of these graphs namely Markov graphs. They correspond to log-linear statistical models of random graphs with general dependence structure and Markov dependence \cite{bremaud} widely used by statisticians and social network analysts \cite{Snijders02}.

To motivate the study of such graphs, we consider the situation where we have measurements of a number of network properties, or observables, for a real-world network, and wish to come up with a network model that exhibits similar properties. Denote these observables by $(x_i)_{i=1,\dots,k}$ and denote by $(\bar{x}_i)_{i=1,\dots,k}$ their measured average value. Let $\mathcal{G}$ a set of graphs, and let $G$ be a graph in $\mathcal{G}$. To describe a family of graphs that reproduce the graph's observed properties, we wish to choose a probability distribution $\mu$ on $\mathcal{G}$ such that
\begin{equation}\label{eq-observables}
\sum_{G\in \mathcal{G}}\mu(G) x_i(G)=\bar{x}_i\:,\quad \forall i=1,\dots,k
\end{equation}
where $x_i(G)$ is the value taken by $x_i$ in the graph $G$. Clearly, there are infinitely many such probability distributions; a popular choice is the one that maximises the Gibbs or Shannon Entropy 
$$S=-\sum_{G\in \mathcal{G}}\mu(G)\log{\mu(G)}\:$$
 subject to (\ref{eq-observables}) and the normalising condition $\sum_{G\in \mathcal{G}}\mu(G)=1$. Introducing Lagrange multipliers one can easily show \cite{Newman04} that the maximum entropy is achieved for the distribution 
\begin{equation}\label{eq-erg}
\mu(G)=\frac{1}{Z}e^{-H(G)},\quad H(G)=\sum_{i=1}^k \theta_i x_i(G)\:,
\end{equation}
and $Z=\sum_{G\in \mathcal{G}}e^{-H(G)}$ is the normalising constant. Graphs drawn according to distributions defined by (\ref{eq-erg}) are called exponential random graphs. Thus, they are random graphs with maximum entropy subject to the specified constraints.

Exponential random graphs can be generated using suitable random walks on the space of graphs, for which they arise as the stationary distribution. More precisely, given $H(G)$, a cost or energy function associated with the graph $G$, define the Markov chain on $\mathcal{G}$ with transition $$p_{G,G'}=\min\left(1,e^{-(H(G')-H(G))}\right)\:.$$ 
It can easily be shown that the transition matrix fulfills the detailed balance condition (the Markov
chain is reversible) and the corresponding stationary distribution is given by the Boltzmann type probability distribution $\mu(G)=Z^{-1}e^{-H(G)}$.

In this paper, we study the particular case of graphs with $n$ nodes and $E$ edges for which $H(G)=\sum_{i=1}d_i^2$, where $d_i$ is the degree of node $i$. This model naturally appears in the context of load balancing in certain communication networks, namely peer-to-peer overlays.  Such overlays are used to support many popular file-sharing applications on the Internet. A primary objective in designing such overlays is to ensure connectivity of the resulting graph even in the face of node and edge disconnections. 
We can model an overlay as a graph with $n$ nodes representing the peers connected by edges describing whether two peers know each other or not. We assume that the ``who knows who'' relationship is symmetric, i.e., the graph is undirected. In \cite{GMK03a}, an algorithm is described that ensures the construction of an Erd\"os-R\'enyi-like overlay, wherein any pair of peers is connected with a given probability independently from other pairs. It is known that such graphs are connected if the mean degree of nodes is of order higher than $\log{n}$ \cite{bollobas}, and the result is true for more general graphs \cite{BB90}. In \cite{GMK03b}, the exponential random graph model with energy function $H(G)=\sum_{i=1}d_i^2$ was proposed as a mechanism for achieving better load balancing and greater resilience to random link failures. \footnote{This paper expands on an earlier short version which appeared in the proceedings of the 41st Allerton Conference on Communications, Control and Computing \cite{GM03}.}
 
The outline of the rest of the paper is as follows. We analyse the degree distribution of such graphs in section \ref{se-degree} and show that the degrees are concentrated around their mean value with high probability ($\whp$). In section \ref{se-cuts}, we derive asymptotic results on the number of edges crossing a graph cut and use these results  $(i)$ to compute the graph expansion and conductance in paragraph \ref{sse-conductance}, and $(ii)$ to analyse the graph resilience to random failures in paragraph \ref{sse-resilience}.

\section{Degree distribution}\label{se-degree}

We work with labelled graphs throughout. We consider the following 
random graph model on $\nv$ nodes with $\ne$ edges:
\begin{equation} \label{eq:graph-distr}
\mu_{\nv}(G) = \frac{1}{Z} \exp{\left(-\beta \sum_{i=1}^{\nv} d_i^2\right)}{\bf 1}_{\{\sum_{i=1}^nd_i=2E\}},
\end{equation}
where $d_i$ denotes the degree of node $i$ in the graph $G$, 
$\beta$ is a specified parameter, and $Z$ is a normalizing 
constant. 

Our aim in this section is to show that graphs generated according to (\ref{eq:graph-distr}), with $2E=cn\log n$, have a sequence of degrees that are concentrated around their mean value. 
 
The probability measure $\mu_{\nv}$ on graphs induces 
a probability measure on degree distributions, which we denote by 
$\pi_{\nv}$. For $\dvec=(d_1,\ldots,d_{\nv})$,
\begin{equation} \label{eq:deg-distr1}
\pi_{\nv}(\dvec) = \frac{1}{Z_{\nv}} G_{\nv}(\dvec) e^{-\beta 
\sum_{i=1}^{\nv} d_i^2} {\bf 1}_{\{\sum_{i=1}^nd_i=2E\}},
\end{equation}
where $G_{\nv}(\dvec)$ is the number of graphs having the degree sequence
$\dvec$, and $Z_{\nv}$ is a normalizing constant. We can rewrite the above as
\begin{eqnarray} 
\pi_{\nv}(\dvec) &=& \frac{1}{Z_{\nv}(\gamma)} \left[ \frac{\ne! 2^{\ne}}
{(2\ne)!} G_{\nv}(\dvec) \prod_{i=1}^{\nv} (d_i!) \right] 
\prod_{i=1}^{\nv} \frac{1}{d_i!} e^{-\beta d_i^2 + \gamma (\log \nv) d_i} 
{\bf 1}_{\{\sum_{i=1}^nd_i=2E\}} \nonumber \\
&=& \frac{\tilde G_{\nv}(\dvec)}{Z_{\nv}(\gamma)} \prod_{i=1}^{\nv}
\frac{1}{d_i!} e^{-\beta d_i^2 + \gamma (\log \nv) d_i}{\bf 1}_{\{\sum_{i=1}^nd_i=2E\}}.  \label{eq:deg-distr2}
\end{eqnarray}
The introduction of the tilt parameter $\gamma$ does not change the
distribution as it multiplies $\pi_{\nv}(\dvec)$ by 
$e^{2 \gamma \ne \log \nv}$. This is a constant since the total number 
of edges is fixed. Thus, it can be absorbed into the normalization factor 
$Z_{\nv}(\gamma)$ along with the term $\ne! 2^{\ne}/(2\ne)!$. 

To construct a graph with a given degree distribution, we use the standard {\em configuration model} \cite{bollobas}: To each node $i$ we associate $d_i$ labelled half-edges, also called configuration points or stubs. All stubs need to be matched to construct the graph, this is done by randomly connecting them. When a stub of $i$ is matched with a stub of $j$, we interpret this as an edge between $i$ and $j$. The graph obtained following this procedure may not be simple, i.e., may contain self-loops due to the matching of two stubs of $i$, and multi-edges due to the existence of more than one matching between two given nodes. 

To restrict ourselves to the family of simple graphs we define the {\em erased configuration model}. Starting from the multigraph obtained through the configuration model, we merge all multiple edges into a single edge and erase all self-loops. It is shown in \cite{vanderHofstad}, that provided that the maximum degree of the graph $d_{\max}$ is such that $d_{\max}=o(\sqrt{n})$, the configuration model and the erased configuration model are asymptotically equivalent, in probability, and every simple graph thus obtained corresponds exactly to $\prod_{i=1}^n d_i!$ distinct configurations describing the number of ways stubs are assigned. We will show in Theorem \ref{thm:main} that the above condition is indeed satisfied.

We denote the minimum and maximum degrees by $d_{\min}$ and $d_{\max}$
respectively. The parameter $\tilde{G}_n({\bf d})$ introduced above corresponds to the probability of obtaining a simple graph in the configuration model. This implies the upper bound 
$\tilde G_{\nv}(\dvec) \le 1$ for any degree sequence $\dvec$. Moreover, if
$d_{\max} = o(\ne^{1/4})$, then McKay and Wormald~\cite{MW90} establish 
the equivalence, for $n$ large,
\begin{equation} \label{eq:graphcount-equiv1}
\tilde G_{\nv}(\dvec) \sim e^{-\lambda-\lambda^2}, \; \hbox{where} \;
\lambda = \frac{1}{4\ne} \sum_{i=1}^{\nv} d_i(d_i-1).
\end{equation}

Given a degree sequence $\dvec$, we define the mean degree $\dmean=
\sum_{i=1}^n d_i/n$ and the variance $\mbox{Var}(\dvec) = \frac{1}{\nv} 
\sum_{i=1}^{\nv} (d_i-\dmean)^2$.
We are interested in a regime where $\dmean = c\log \nv$
for some specified constant $c$, so that $\ne = c\nv \log \nv/2$.

For fixed constants $\alpha_1$ and $\alpha_2$, we define 
the following sets of degree sequences:
\begin{eqnarray*}
A &=& \{ \dvec : \; \dmean = c\log \nv \}, \\
A_1(\alpha_1,\alpha_2) &=& \{ \dvec : \; - \sqrt{\alpha_1 \log \nv} \le 
d_i - \dmean \le \sqrt{\alpha_2 \log \nv} \:, \forall i=1,\dots,\nv \}, %\\
\end{eqnarray*}
Note that, in the regime $\dmean = c\log \nv$, $\pi_{\nv}$ is supported on $A$,
and so $\pi_{\nv}(B)=\pi_{\nv}(A\cap B)$ for any set $B$ of labelled
graphs on $\nv$ nodes. Define $\hat A_1(\alpha_1,\alpha_2) = A \cap A_1(\alpha_1,\alpha_2)$. 
We wish to show that 

\begin{theorem} \label{thm:main}
There exist constants $\alpha_1, \alpha_2$ such that 
$\pi_n(\hat A_1(\alpha_1,\alpha_2))$ goes to 1 as $\nv$ goes to infinity.
\end{theorem}

The above theorem states that for the random graph model defined by the distribution (\ref{eq:graph-distr}), the node degrees concentrate about their mean value. Specifically,
 all node degrees are within order 
$\sqrt{\log \nv}$ of the mean,  $\whp$. This is in contrast to the 
Erd\"os-R\'enyi model (with the same number of edges) where the maximum fluctuation of node degrees is typically of order 
$\log \nv$. The rest of the section is devoted to the proof of this theorem. To this end, we start by proving that 

\begin{theorem} \label{thm:estimate}
Define the event $A_2 = \{ \dvec : \; d_i \le \nv^{1/4} \:, \forall i=1,\dots,\nv \}$. Then
$$\pi_n(A_2^c)\to 0,\qquad \text{as }\nv \to \infty\:,$$
and the estimate in (\ref{eq:graphcount-equiv1}) holds. \end{theorem}

To prove this we first state a series of lemmas which are proved in Appendix \ref{se-appendix}.
 
If $\dvec \in
\hat A_1(\alpha_1,\alpha_2)$, then in the regime $\ne = c\nv \log \nv/2$, we have $d_{\max} = o(\ne^{1/4})$. Observe from (\ref{eq:graphcount-equiv1}) that $4\ne \lambda = 
\nv(\mbox{Var}(\dvec)+\dmean^2-\dmean)$. Moreover, for $\dvec \in
\hat A_1(\alpha_1,\alpha_2)$, we have  $\mbox{Var}(\dvec) \le \max \{ \alpha_1,
\alpha_2 \} \log \nv$, so that 

\begin{equation}\label{eq-lambdaupper}
\lambda \le \frac{1}{2}( c\log \nv - 1 
+ \frac{1}{c} \max\{ \alpha_1,\alpha_2 \} )\:.
\end{equation}
Hence, 
\begin{equation}\label{eq-G_nequi}
\dvec \in \hat A_1(\alpha_1,\alpha_2) \; \Rightarrow \; 
\frac{1}{\tilde G_{\nv}(\dvec)} 
\sim e^{\lambda+\lambda^2} \le e^{\frac{c^2 \log^2 \nv}{2}},
\end{equation}
for all $\nv$ sufficiently large. Recall that $\tilde G_{\nv}(\dvec)
\le 1$ for all $\dvec$ and, in particular, for $\dvec \in A_2^c$, the
complement of $A_2$. Thus, it follows from (\ref{eq:deg-distr2}) and (\ref{eq-G_nequi}) that, 
for $\nv$ sufficiently large,
\begin{equation} \label{eq:tail2}
\frac{\pi_n(A_2^c)}{\pi_n(\hat A_1(\alpha_1,\alpha_2))} \le 
e^{ \frac{c^2 \log^2 \nv}{2} }
\frac{ \sum_{\dvec \in A_2^c} \prod_{i=1}^{\nv} \frac{1}{d_i!}
e^{-\beta d_i^2+\gamma (\log \nv) d_i} }{ \sum_{\dvec \in 
\hat A_1(\alpha_1,\alpha_2)} \prod_{i=1}^{\nv} \frac{1}{d_i!} 
e^{-\beta d_i^2+\gamma (\log \nv) d_i} }.
\end{equation}

Let $D_1,\ldots,D_{\nv}$ be independent and identically distributed (iid)
random variables, with
\begin{equation} \label{eq:marginal}
\PP(D_1 = k) = \frac{1}{F(\gamma)} \frac{1}{k!} 
e^{-\beta k^2+\gamma (\log \nv) k},
\quad k\in \Nats,
\end{equation}
where $F(\gamma)$ is a normalization constant. The dependence of the $D_i$
on $\nv$ and $\gamma$ has not been made explicit in the notation. We choose
$\gamma$ so that $\EE D_1 = c\log \nv$, for a specified constant, $c$; this is
possible by the following lemma. 

\begin{lemma}\label{le-deg dist1}
Let $x_{\gamma}=\frac{1}{2\beta}\left(\gamma \log{n}+\log\log{n}+\frac{\gamma}{2\beta}\right)$, and let $k_{\gamma}-1$ denote the integer part of $x_{\gamma}$. Then, $\EE D_1-k_{\gamma}$ and $\mbox{Var}(D_1)$ remain bounded as $n$ tends to $\infty$. Moreover, let $\alpha=2\beta\left(x_{\gamma}-k_{\gamma}+\frac{1}{2}\right)$ and 
$$\psi(\theta)=\frac{\sum_{j=-\infty}^{\infty} e^{\theta j- \beta j^2} }{\sum_{j=-\infty}^{\infty}e^{- \beta j^2} }\:.$$

Then, the moment generating function of $D_1$ satisfies
$$\EE\left[e^{\theta D_1}\right]\sim e^{\theta k_{\gamma}} \frac{\psi(\theta+\alpha)}{\psi(\alpha)},\qquad \text{as }\nv \to \infty\\:.$$
\end{lemma}

\noindent \emph{Proof}: 
See proof in section \ref{proof-deg dist1}.
\hfill $\Box$

Let ${\bf D}$ denote the random vector $(D_1,\ldots,D_{\nv})$.
We can now rewrite (\ref{eq:tail2}) as
\begin{equation} \label{eq:sloppybd1}
\frac{\pi_n(A_2^c)}{\pi_n(\hat A_1(\alpha_1,\alpha_2))} \le 
e^{ \frac{c^2 \log^2 \nv}{2} }
\frac{ \PP( {\bf D} \in A_2^c) }{ \PP( {\bf D} \in
\hat A_1(\alpha_1,\alpha_2)) }.
\end{equation}

\begin{lemma} \label{lem:tail1}
There exists a constant $K>0$, independent of $\nv$, such that 
\begin{equation}\label{eq-A_2}
\PP({\bf D} \in A_2^c) \le 
Kne^{-\beta \sqrt{\nv}/4}\:.
\end{equation}
\end{lemma}

\noindent \emph{Proof}: 
See proof in section \ref{proof-le-tail1}.
\hfill $\Box$

Let $(\tilde D_1,\ldots,\tilde D_{\nv})$ have the joint distribution
of $(D_1,\ldots,D_{\nv})$ conditional on ${\bf D} \in
A_1(\alpha_1,\alpha_2)$. Equivalently, $\tilde D_1,\ldots,\tilde D_{\nv}$
are iid, with $\tilde D_j$ having the distribution of $D_j$ conditional on
$$- \sqrt{\alpha_1 \log \nv} \le D_j-\EE D_j \le  \sqrt{\alpha_2\log \nv}\:.$$ Now
\begin{eqnarray}
\PP\left({\bf D} \in \hat A_1(\alpha_1,\alpha_2)\right) 
&=&
\PP\left({\bf D} \in A_1(\alpha_1,\alpha_2)\right)~\PP\left(\sum_{j=1}^{\nv} D_j = 
c\nv \log \nv \mid {\bf D} \in A_1(\alpha_1,\alpha_2)\right) \nonumber \\
& =& \; \PP\left({\bf D} \in A_1(\alpha_1,\alpha_2)\right)~
\PP\left(\sum_{j=1}^{\nv} \tilde D_j = c\nv \log \nv\right) \label{eq:denom-sloppy1}
\end{eqnarray}
Suppose $\alpha_1, \alpha_2 > 0$ are chosen large enough so that, for $n$ large, $\EE\tilde D_1 = \EE D_1
= c\log \nv$. 

We wish to estimate the probability that $\tilde{D_1}+\tilde{D_2}+\cdots+\tilde{D_{\nv}} = 
c\nv \log \nv$. We shall do this using a result from~\cite{M79}.
For $j=1,\ldots,\nv$, define the centred random variables,
$X_{\nv j}=\tilde{D_j} -\EE \tilde{D_j}$; we have made the dependence of the distribution 
of $\tilde{D_j}$ on $\nv$ explicit in the notation.
Thus, $X_{\nv 1}, X_{\nv 2},\ldots, X_{\nv \nv}$ is an array of
integer-valued zero mean random variables such that, for each $\nv$,
$X_{\nv 1},\ldots,X_{\nv \nv}$ are independent and identically distributed.
Now, to apply~\cite[Theorem 1]{M79}, we need the following result.

\begin{lemma} \label{lem:cgf}
The random variables, $\{ X_{\nv j}, j=1,\ldots,\nv, \nv \in \Nats \}$,
satisfy the following conditions:\\
(i) $\limsup_{\nv \to \infty} \EE[e^{\theta |X_{\nv 1}|}] < \infty$ for
some $\theta > 0$. \\
(ii) $\liminf_{\nv \to \infty} \mbox{Var}(X_{\nv 1}) > 0$. \\
(iii) $\liminf_{\nv \to \infty} \sum_{j=-\infty}^{\infty}
\min \{ \PP(X_{\nv 1}=j), \PP(X_{\nv 1}=j+1) \} > 0$.
\end{lemma}

\noindent \emph{Proof}: 
See proof in section \ref{proof-le-cgf}.
\hfill $\Box$

Indeed, an immediate corollary of \cite[Theorem 1]{M79} is that
\begin{theorem}
If a sequence of independent random variables, $\{ X_{\nv j}, j=1,\ldots,\nv, \nv \in \Nats \}$,
satisfies the conditions (i),  (ii) and (iii)  of Lemma \ref{lem:cgf}, then
\begin{equation*} 
\PP\left(\sum_{j=1}^{\nv}   X_{\nv j}= \sum_{j=1}^{\nv}\EE X_{\nv j}\right) = 
\frac{1}{\sqrt{2\pi\sum_{j=1}^{\nv} \mbox{Var}(X_{\nv j})}} \Bigl( 1+ O\Bigl( \frac{1}{\nv} 
\Bigr) \Bigr),
\end{equation*}
\end{theorem}

A direct application of the above result yields
\begin{equation} \label{eq:locallimit1}
\PP\left(\sum_{j=1}^{\nv} \tilde D_j = c\nv \log \nv\right) = 
\frac{1}{\sqrt{2\pi\nv} \tilde \sigma} \Bigl( 1+ O\Bigl( \frac{1}{\nv} 
\Bigr) \Bigr),
\end{equation}
where $\tilde \sigma = \mbox{Var}(\tilde D_1)$ remains bounded as
$\nv \to \infty$. Combining this with (\ref{eq:sloppybd1}), (\ref{eq-A_2}) and
(\ref{eq:denom-sloppy1}), we get
\begin{equation} \label{eq:sloppybd2}
\pi_n(A_2^c) \le
\frac{\pi_n(A_2^c)}{\pi_n(\hat A_1(\alpha_1,\alpha_2))} \le
e^{ \frac{c^2 \log^2 \nv}{2} }
\frac{ \sqrt{2\pi} \tilde \sigma K \nv^{3/2} e^{-\beta \sqrt{\nv}/4} }
{ \PP({\bf D} \in A_1(\alpha_1,\alpha_2)) } \Bigl( 1+ O\Bigl( \frac{1}{\nv} 
\Bigr) \Bigr).
\end{equation}

\begin{lemma} \label{lem:a1prob}
Let ${\bf D}$ denote the random vector $(D_1,\ldots,D_{\nv})$. Given
any $K>0$, we can choose $\alpha_1$ and $\alpha_2$ such that
$\PP({\bf D} \in A_1(\alpha_1,\alpha_2)^c) < e^{-K\log \nv}$ for all
$\nv$ sufficiently large.
\end{lemma}

\noindent \emph{Proof}: 
See proof in section \ref{proof:a1prob}.
\hfill $\Box$

Combining the above lemma with the bound in (\ref{eq:sloppybd2}), it is immediate that 
$\pi_n(A_2^c)\to 0$ as $\nv \to \infty$ which establishes the claim of Theorem \ref{thm:estimate}.
Thus, to prove Theorem \ref{thm:main}, we can restrict our attention to graphs with degree sequences in $A_2$, for which we can use the estimate in (\ref{eq:graphcount-equiv1}).

\noindent \emph{Proof of Theorem \ref{thm:main}}:
Observe that
\begin{eqnarray*}
\pi_n(\hat A_1(\alpha_1,\alpha_2)) 
&=& \pi_{\nv}(A) -\pi_{\nv}(A\setminus A_1(\alpha_1,\alpha_2)) \\
&\ge& \pi_{\nv}(A) -\pi_{\nv}((A\setminus A_1(\alpha_1,\alpha_2))\cap A_2)
-\pi_{\nv}(A_2^c).
\end{eqnarray*}
But $\pi_{\nv}(A)=1$ by definition, and we have shown above that
$\pi_n(\hat A_2^c) \to 0$ as $\nv \to \infty$. Hence, it suffices to 
show that 
\begin{equation} \label{eq:suff}
\pi_n((A\setminus A_1(\alpha_1,\alpha_2))\cap A_2) \to 0, \quad
\mbox{ as $\nv\to \infty$.} 
\end{equation}
Recall from (\ref{eq:graphcount-equiv1}) that, if 
$\dvec \in A_2$, then $\tilde G_{\nv}(\dvec) \sim 
e^{-\lambda(\dvec)-\lambda(\dvec)^2}$. Now, 
\[
\lambda(\dvec) = \frac{\mbox{Var}(\dvec)+{\overline d}^2-{\overline d}}
{2\overline d} \ge \frac{c\log \nv-1}{2}, \quad \forall \: \dvec \in A,
\]
since the mean degree, ${\overline d}=c\log \nv$. In particular,
the above lower bound on $\lambda(\dvec)$ holds for all degree sequence $\dvec$ in
$(A\setminus A_1(\alpha_1,\alpha_2)) \cap A_2$, since this is a 
subset of $A$. 

In addition, we saw earlier  in (\ref{eq-lambdaupper}) that, if $\dvec \in 
\hat A_1(\alpha_1,\alpha_2)$, then
\[
\lambda(\dvec) \le \frac{1}{2} \Bigl( c\log \nv - 1 + \frac{1}{c}\max
\{ \alpha_1,\alpha_2 \} \Bigr)\:,
\]

and the estimate in (\ref{eq:graphcount-equiv1}) holds.

Denote $\max \{ \alpha_1,\alpha_2 \}$ by $\alpha$. 
Now, by (\ref{eq:deg-distr2}),
\begin{eqnarray*}
\frac{ \pi_n\left((A\setminus A_1(\alpha_1,\alpha_2)) \cap A_2\right) }
{ \pi_n(\hat A_1(\alpha_1,\alpha_2)) } 
&=& \; \frac{ \sum_{ \dvec \in (A\setminus A_1(\alpha_1,\alpha_2)) \cap A_2 }
e^{ -\lambda(\dvec)-\lambda(\dvec)^2 } \prod_{i=1}^{\nv}
\frac{1}{d_i!} e^{ -\beta d_i^2+\gamma (\log \nv) d_i } }
{ \sum_{ \dvec \in \hat A_1(\alpha_1,\alpha_2) } 
e^{ -\lambda(\dvec)-\lambda(\dvec)^2 } \prod_{i=1}^{\nv}
\frac{1}{d_i!} e^{ -\beta d_i^2+\gamma (\log \nv) d_i } } \\
&\leq & \; e^{ \frac{\alpha}{2c} \left(
c\log \nv + \frac{ \alpha }{2c} \right) }
\frac{ \sum_{ \dvec \in (A\setminus A_1(\alpha_1,\alpha_2)) \cap A_2 }
\frac{1}{d_i!} e^{ -\beta d_i^2+\gamma (\log \nv) d_i } }
{ \sum_{ \dvec \in \hat A_1(\alpha_1,\alpha_2) }
\frac{1}{d_i!} e^{ -\beta d_i^2+\gamma (\log \nv) d_i } }.
\end{eqnarray*}
In other words, there are constants $\kappa_1$ and $\kappa_2$ such that
\begin{eqnarray}
\frac{ \pi_n((A\setminus A_1(\alpha_1,\alpha_2)) \cap A_2) }
{ \pi_n(\hat A_1(\alpha_1,\alpha_2)) } 
&\le& \kappa_1 e^{\kappa_2 \log \nv} 
\frac{ \PP( {\bf D} \in (A\setminus A_1(\alpha_1,\alpha_2)) \cap A_2 ) }
{  \PP( {\bf D} \in \hat A_1(\alpha_1,\alpha_2) ) } \nonumber \\
&\le& \kappa_1 e^{\kappa_2 \log \nv}
\frac{ \PP( {\bf D} \in A\setminus A_1(\alpha_1,\alpha_2) ) }
{  \PP( {\bf D} \in \hat A_1(\alpha_1,\alpha_2) ) }. \label{eq:main-bd1}
\end{eqnarray}
Now, by Lemma~\ref{lem:a1prob}, for any given $K>0$, we can choose
$\alpha_1$ and $\alpha_2$ such that $\PP\left( {\bf D} \in 
A_1(\alpha_1,\alpha_2)^c \right) \le e^{-K\log \nv}$. Thus,
\begin{equation} \label{eq:main-num}
\PP({\bf D} \in A\setminus A_1(\alpha_1,\alpha_2)) \le
\PP({\bf D} \in A_1(\alpha_1,\alpha_2)^c) \le e^{-K\log \nv}.
\end{equation}
Moreover, analogous to (\ref{eq:locallimit1}), we have
\[
\PP({\bf D} \in A) = \PP(\sum_{j=1}^{\nv} D_j = c\nv \log \nv)
= \frac{1}{ \sqrt{2\pi\nv} \sigma } \left( 1+O\left( \frac{1}{n}
\right) \right),
\]
where $\sigma = \mbox{Var}(D_1)$ remains bounded as $\nv \to \infty$.
Therefore,
\begin{eqnarray} 
\PP({\bf D} \in \hat A_1(\alpha_1,\alpha_2))
&=& \PP({\bf D}\in A) - \PP({\bf D} \in A\cap A_1(\alpha_1,\alpha_2)^c) 
\nonumber \\
&\ge& \PP({\bf D}\in A) - \PP({\bf D} \in A_1(\alpha_1,\alpha_2)^c) \nonumber \\
&=& \frac{1}{ \sqrt{2\pi\nv} \sigma } \left( 1+O\left( \frac{1}{n}
\right) \right).  \label{eq:main-den}
\end{eqnarray}
Substituting (\ref{eq:main-num}) and (\ref{eq:main-den}) in
(\ref{eq:main-bd1}), we have
\begin{eqnarray*}
\pi_n((A\setminus A_1(\alpha_1,\alpha_2)) \cap A_2) 
&\le& \frac{ \pi_n(A\setminus A_1(\alpha_1,\alpha_2)) }
{ \pi_n(\hat A_1(\alpha_1,\alpha_2)) } \\
&\le& \kappa_1 \sigma \sqrt{2\pi \nv} \, e^{(\kappa_2-K)\log \nv} 
\left( 1+O\left( \frac{1}{n} \right) \right).
\end{eqnarray*}
Since $K$ can be chosen arbitrarily large, the above quantity
goes to zero as $\nv \to \infty$, which establishes (\ref{eq:suff})
and the claim of the theorem. \hfill $\Box$

\section{Graph cuts}\label{se-cuts}

Given a graph $G$ and a subset $U$ of its vertex set, let $e_U(G)$ denote
the number of edges incident within $U$ (i.e., having both their vertices
with $U$); let $e_{U,U^c}(G)$ denote the number of edges having one
vertex in $U$ and the other in its complement, $U^c$ (i.e., crossing the cut $(U,U^c)$); and denote by $u$ or $|U|$ the number of vertices or size of $U$.

Let ${\bf d}(G) = (d_1,d_2,\ldots,d_{\nv})$ denote the degree sequence of $G$ and define the volume of a subset of vertices $U$ by
$$\mbox{Vol}(U) = \sum_{i\in U} d_i.$$ Note that 
\begin{equation}\label{eq-vol}
2e_U(G) + e_{U,U^c}(G) = \mbox{Vol}(U).
\end{equation}
In the remainder of this section we derive lower bounds for the graph cuts. To this end we will show that there exists a constant $\delta$ such that $e_{U,U^c}(G)>(1-\delta)|U|c\log n,\:\whp$, using different techniques depending on the size of $U$, when $|U|\leq n/2$.

\begin{prop}\label{expander0}
For any $\epsilon>0$, there exists $\delta_1\in(0,1)$, independent of $n$, such that, if the subset of vertices $U$ is such that $u\leq 2\epsilon c\log \nv$, then $e_{U,U^c}(G) \ge (1-\delta_1) u c\log \nv$, $\whp$.
\end{prop}
\noindent \emph{Proof}: 
Denote $|U|$ by $u$.
Suppose first that $u \le 2\epsilon c\log \nv$, for a given
$\epsilon>0$. The number of edges incident within $U$ can be at most 
${u\choose 2}$, so $e_U(G) \le \epsilon u c\log \nv$, for all $U$.
Now, for any degree sequence ${\bf d} \in A_1(\alpha_1,\alpha_2)$,
$ \mbox{Vol}(U) \ge cu\log \nv - u\sqrt{\alpha_1 \log \nv}$. By Theorem \ref{thm:main}, it is not restrictive to consider only graphs with degree sequences belonging to the set $A_1(\alpha_1,\alpha_2)$. Hence, using (\ref{eq-vol}) for graphs
$G$ with such degree sequences,
\[
e_{U,U^c}(G) \ge u[(1-2\epsilon) c\log \nv - \sqrt{\alpha_1 \log \nv}].
\]
Let $\delta_1=3\epsilon$. Then, for $\nv$ sufficiently large,
$e_{U,U^c}(G) \ge (1-\delta_1) u c\log \nv$, $\whp$, whenever 
$u \le 2\epsilon c\log \nv$ and the claim of the proposition is established.  
\hfill $\Box$

To prove a similar result for all subsets $U$ such that $u\leq n/2$ we will use the configuration model \cite{bollobas}. Fix a degree sequence ${\bf d} \in A_1(\alpha_1,\alpha_2)$.
By (\ref{eq:graph-distr}), all graphs with the same degree sequence 
are equally likely under the distribution $\mu_{\nv}$, so we can use 
the configuration model to generate a random graph with 
this distribution, conditional on the degree sequence. 

%Recall that each graph $G$ with a given degree sequence corresponds to
%the same number of configurations, namely $\prod_{i=1}^{\nv} d_i!$,
%but a configuration may not yield a simple graph (it could have loops 
%and multiple edges). Assuming that every ``bad" configuration (namely,
%a configuration $H$ with $e_{U,U^c}(H) < (1-\delta)uc\log \nv$ for
%some $U\subseteq V$) corresponds to a simple graph yields an upper bound 
%on the fraction of graphs $G$ which are bad. 

For constants $\delta\in(0,1)$, $\epsilon>0$, and $\tau>0$, for $\nv \in \Nats$
and a degree sequence ${\bf d}$,
we define the following subsets of graphs on a vertex set $V$ of
cardinality $\nv$:
\begin{eqnarray}
{\cal E}_1(\nv,\delta,\tau,{\bf d}) &=&
\{ G: {\bf d}(G)={\bf d} \mbox{ and } e_{U,U^c}(G) < (1-\delta) u c\log \nv 
\nonumber \\
&& \quad \mbox{for some $U\subseteq V$ with $2\epsilon c\log \nv < u \le \tau \nv$} \}, 
\label{badgraphs1} \\
{\cal E}_2(\nv,\delta,\tau,{\bf d}) &=&
\{ G: {\bf d}(G)={\bf d} \mbox{ and } e_{U,U^c}(G) < (1-\delta) u c\log \nv  \nonumber \\
&& \quad \mbox{for some $U\subseteq V$ with $\tau \nv < u \le \nv/2$} \}. 
\label{badgraphs2}
\end{eqnarray}
We also define
\begin{equation} \label{badgraphs}
{\cal E}_1(\nv,\delta,\tau) = 
\bigcup_{\bf d} {\cal E}_1(\nv,\delta,\tau,{\bf d}), \quad
{\cal E}_2(\nv,\delta,\tau) =
\bigcup_{\bf d} {\cal E}_2(\nv,\delta,\tau,{\bf d}).
\end{equation}
We shall derive bounds on the probabilities of these sets using the
configuration model~\cite{bollobas}. To this end, we define the analogous
sets of configurations $\hat {\cal E}_1(\nv,\delta,\tau,{\bf d})$,
$\hat {\cal E}_2(\nv,\delta,\tau,{\bf d})$, $\hat {\cal E}_1(\nv,\delta,\tau)$
and $\hat {\cal E}_2(\nv,\delta,\tau)$. It is useful to define the
following sets of configurations on the same vertex set. More precisely, given a degree
sequence ${\bf d}=(d_1,d_2,\ldots,d_{\nv})$, and for $H$ a configuration on $V$, we define
\begin{eqnarray}
\hat{\cal E}_1(\nv,\delta,\tau,{\bf d}) &=& 
\{ H: {\bf d}(H)={\bf d} \text{ and }
e_{U,U^c}(H) < (1-\delta) u c\log \nv \nonumber \\ &&\quad \mbox{ for some } U\subseteq V \mbox{ with } 
2\epsilon c\log \nv < u \le \tau \nv \}, \label{badconfigs1}\\
\hat{\cal E}_2(\nv,\delta,\tau,{\bf d}) &=& 
\{ H: {\bf d}(H)={\bf d} \text{ and }
e_{U,U^c}(H) < (1-\delta) u c\log \nv \nonumber \\ &&\quad \mbox{ for some } U\subseteq V \mbox{ with } 
\tau \nv < u \le \nv/2 \}, \label{badconfigs2}
\end{eqnarray}
and
\begin{equation} \label{badconfigs}
\hat{\cal E}_1(\nv,\delta,\tau) = 
\bigcup_{\bf d} \hat{\cal E}_1(\nv,\delta,\tau,{\bf d}), \quad
\hat{\cal E}_2(\nv,\delta,\tau) =
\bigcup_{\bf d} \hat{\cal E}_2(\nv,\delta,\tau,{\bf d}).
\end{equation}

Recall that configurations correspond to multigraphs, i.e, there may
be loops or multiple edges. A multiple edge is counted the corresponding
number of times in the above definitions.

Since ${\bf d} \in A(\alpha_1,\alpha_2)$, estimate (\ref{eq:graphcount-equiv1}) holds. Using the enumeration formula 
of McKay and Wormald~\cite{MW90}, this bound says that, for $i=1,2$
\begin{equation} \label{configbound1}
\mu_{\nv}({\cal E}_i(\nv,\delta,\tau,{\bf d})|{\bf d}) \le
e^{\lambda+\lambda^2}
\P(H \in \hat {\cal E}_i(\nv,\delta,\tau,{\bf d})|{\bf d}),
\end{equation}
where $\P(\cdot|{\bf d})$ denotes the probability with respect to the 
uniform distribution on configurations with degree sequence ${\bf d}$.
Recall that $\lambda$ was defined in (\ref{eq:graphcount-equiv1})
to be $\sum_{i=1}^{\nv} d_i(d_i-1)/4E$, where $E$ is the number of
edges, i.e., $2E=\sum_{i=1}^{\nv} d_i$. The dependence of $\lambda$
on ${\bf d}$ has been suppressed for notational convenience.
\begin{prop} \label{expander1}
If $\tau\in(0,\frac{1}{1+4e})$, then there exists $\delta_2\in(0,1)$, independent of $n$, such that 
\[
\lim_{\nv \to \infty} \mu_{\nv} ({\cal E}_1(\nv,\delta_2,\tau)) = 0,
\]
where the distribution $\mu_{\nv}$ was defined in (\ref{eq:graph-distr}). 
\end{prop}

\noindent \emph{Proof}:

For degree sequences ${\bf d} \in A_1(\alpha_1,\alpha_2)$ and
any subset $U$ of the vertex set, $\mbox{Vol}(U)\sim uc\log \nv$, for $n$ large.

Hence, by (\ref{eq-vol}), 
$e_{U,U^c}(H) <u(1-\delta) c\log \nv$ for a subset $U$ implies that
$e_U(H) > \frac{\delta}{2}  \mbox{Vol}(U)$, for sufficiently large $\nv$. To prove the proposition it therefore sufffices to show that there exists $\delta_2 \in(0,1)$ such that $\P(e_U(H) > \frac{\delta_2}{2} ~\mbox{Vol}(U))$ tends to $0$ when $n$ tends to infinity.

Recall that for subset $U$ of $V$, the volume of $U$ is given by $\mbox{Vol}(U)= \sum_{i\in U} d_i$. As the half-edges in the configuration model are matched uniformly, $e_U(H)$, the number of
edges incident within $U$ in a random configuration, is bounded above
by a binomial random variable $X$ with parameters $\mbox{Vol}(U)$ and
$\mbox{Vol}(U)/(2E-\mbox{Vol}(U))$. The dependence of $X$ on $U$ has been suppressed
for notational convenience. For $\delta\in(0,1)$, by Chernoff's bound,  we have
\begin{eqnarray*} 
\log \P\left(X > \frac{\delta}{2}~\mbox{Vol}(U)\right) 
& \le &-\mbox{Vol}(U) \Bigl[ \frac{\delta}{2} \log \frac{\frac{\delta}{2} (2E-\mbox{Vol}(U))}{\mbox{Vol}(U)} 
+ (1-\frac{\delta}{2}) \log \frac{(1-\frac{\delta}{2})(2E-\mbox{Vol}(U))}{2E-2\mbox{Vol}(U)} \Bigr] \\
& \le & -\mbox{Vol}(U) \Bigl[ \frac{\delta}{2} \log \frac{\frac{\delta}{2} (2E-\mbox{Vol}(U))}{\mbox{Vol}(U)}
+ (1-\frac{\delta}{2}) \log (1-\frac{\delta}{2}) \Bigr]. 
\end{eqnarray*}

Applying the inequality $\log x \le x-1$ for $x\geq 1$ to $x=1/(1-\frac{\delta}{2})$, we have that
$\log(1-\frac{\delta}{2}) \ge -\frac{\delta}{2}/(1-\frac{\delta}{2})$.
Using the fact that
$$\Bigm| \frac{ \mbox{Vol}(U)}{uc\log \nv}-1 \Bigm| < \frac{\sqrt{\alpha}}{c} 
\frac{1}{\sqrt{\log \nv}}\:,$$
we have
\begin{equation} \label{binom-ubd1}
\log \P\left(X > \frac{\delta}{2} \mbox{Vol}(U)\right) \le -u c\log \nv 
\left[ \frac{\delta}{2} \log\left(\frac{\delta(\nv-u)}{2u}\right) -\frac{\delta}{2} \right]
\Bigl( 1+O \Bigl( \frac{1}{\sqrt{\log \nv}} \Bigr) \Bigr).
\end{equation}

Suppose first that $2\epsilon c\log \nv < u\le \sqrt{n}$. 

For all $\nv$ sufficiently large, equation (\ref{binom-ubd1}) becomes
\[
\log {\P}\left(X > \frac{\delta}{2} \mbox{Vol}(U)\right) \le -\frac{u \delta c}{6} \log^2 \nv.
\]

Since $X$ stochastically dominates $e_U(H)$ (conditional on ${\bf d}$),
we have by the union bound that, for $\nv$ sufficiently large,
\begin{eqnarray}
\P\left(\exists U,\: 2\epsilon c\log \nv < u
\leq \sqrt{\nv},\: e_U(H) > \frac{\delta}{2} \mbox{Vol}(U)\right) \nonumber 
& \le &\sum_{u=2\epsilon c\log \nv}^{\sqrt{\nv}} {\nv \choose u} 
\exp \Bigl( -\frac{u \delta c}{6} \log^2 \nv \Bigr) \nonumber \\
& \le&\sum_{u=2\epsilon c\log \nv}^{\sqrt{\nv}} \frac{1}{u!}
\exp \Bigl( u\log \nv -\frac{u \delta c}{6} \log^2 \nv \Bigr) \nonumber \\
& \le &\kappa_3 \exp \Bigl( -\kappa_4 \epsilon \delta c^2 \log^3 \nv \Bigr)\:,
\label{binom-ubd2}
\end{eqnarray}
for two constants $\kappa_3, \kappa_4>0$. We have used the inequality ${\nv \choose u} \le \nv^u/u!$ to obtain the
second inequality above.

Next, consider $\sqrt{\nv} < u \le \tau \nv$. 

In this case equation (\ref{binom-ubd1}) becomes,

\begin{equation*} 
\log \P\left(X > \frac{\delta}{2} \mbox{Vol}(U)\right) \le -\frac{1}{2}u c\log \nv 
\left[ \delta \log\left(\frac{\delta(1-\tau)}{2\tau}\right) -\delta \right]
\Bigl( 1+O \Bigl( \frac{1}{\sqrt{\log \nv}} \Bigr) \Bigr).
\end{equation*}

If $\tau<\frac{1}{1+4e}$, then there exists $\delta_2\in(0,1)$ such that
\[
\delta_2\log\left( \frac{\delta_2(1-\tau)}{2\tau}\right) -\delta_2 > \frac{2}{c}.
\]
and subsequently, for all $\nv$ sufficiently large 
and for $u\le \tau \nv$, we have 
\[
\log \P\left(X > \frac{\delta_2}{2} \mbox{Vol}(U)\right) \le -2u\log \nv.
\]
Hence, by the union bound,
\begin{eqnarray}
\P\left(\exists U: \sqrt{\nv} < u
< \tau \nv,\: e_U(H) > \frac{\delta_2}{2} \mbox{Vol}(U) \right) \nonumber 
& \leq & \sum_{u=\sqrt{\nv}}^{\tau \nv} {n\choose u} e^{-2u\log \nv}\\ \nonumber 
&\leq & \sum_{u=\sqrt{\nv}}^{\tau \nv} \frac{1}{u!} e^{-u\log \nv}\\
&\leq &   \kappa_5 e^{-\sqrt{\nv}\log n}.  \label{binom-ubd3} 
\end{eqnarray}

By (\ref{configbound1}), (\ref{binom-ubd2}) 
and (\ref{binom-ubd3}), for $n$ large, we can find two constants $\kappa_6, \kappa_7>0$ such that
\[
\mu_{\nv}({\cal E}_1(\nv,\delta_2,\tau,{\bf d})|{\bf d}) \le
e^{\lambda+\lambda^2}\kappa_6 e^{-\kappa_7 \log^3\nv}.
\]
Since $\lambda=O(\log \nv)$,
it is readily checked that $\mu_{\nv}({\cal E}_1(\nv,\delta_2,\tau,{\bf d})|{\bf d} \in
A_1(\alpha_1,\alpha_2) )$ goes to $0$ as $\nv \to \infty$. 

By Theorem~\ref{thm:main},
$\mu_{\nv}( {\bf d} \notin A_1(\alpha_1,\alpha_2) )$ goes to 0 as well. Noting that
\[
\mu_{\nv}({\cal E}_1(\nv,\delta,\tau)) \le
\mu_{\nv}({\cal E}_1(\nv,\delta,\tau,{\bf d})|{\bf d} \in
A_1(\alpha_1,\alpha_2) ) + \mu_{\nv}( {\bf d} \notin A_1(\alpha_1,\alpha_2) ),
\]
the claim of the proposition is established.  \hfill $\Box$

 Next, we find a similar lower bound for $e_{U,U^c}(G)$ that
holds, $\whp$, for subsets $U$ with $\tau \nv < u \le \nv/2$.

\begin{prop} \label{expander2}
For $\tau>0$, there exists $\delta_3\in(0,1)$, independent of $n$, such that 
\[
\lim_{\nv \to \infty} \mu_{\nv} ({\cal E}_2(\nv,\delta_3,\tau)) = 0.
\]
\end{prop}

\noindent\emph{Proof}:
As in the proof of Proposition~\ref{expander1}, we fix a degree sequence ${\bf d}$
and a subset $U$, and bound the probability that $e_{U,U^c}(G) < u(1-\delta)c\log n$ in
terms of the probability that $e_{U,U^c}(H) < u(1-\delta)c\log n$, where $H$ is drawn
uniformly at random from configurations with degree sequence ${\bf d}$, i.e.,
\begin{equation} \label{configbound2}
\mu_{\nv}({\cal E}_2(\nv,\delta,\tau,{\bf d})|{\bf d}) \le
e^{\lambda+\lambda^2}
\P(H \in \hat {\cal E}_2(\nv,\delta,\tau,{\bf d})|{\bf d}),
\end{equation}

Fix constants $\tau>0$ and $\delta\in(0,1)$, and a degree sequence ${\bf d}$. 
Let $U$ be a subset of the vertex set with 
$\tau \nv < u \le \nv/2$, and let $j<(1-\delta)uc\log n\leq \frac{1}{2}(1-\delta) cn \log n$. Recall that the number of 
configurations with degree sequence ${\bf d}$ is
\begin{equation} \label{nconfig}
H_{\nv}({\bf d}) = \frac{(2E)!}{E!2^E} \prod_{i=1}^{\nv} d_i!,
\end{equation}
where $E = \sum_{i=1}^{\nv} d_i/2$ is the total number of edges. The
number of these configurations with exactly $j$ edges crossing the
cut between $U$ and $U^c$ is
{\small
\begin{equation} \label{nconfigcut}
{\cal H}_{U,U^c}(j) 
\le {\mbox{Vol}(U) \choose j}{2E-\mbox{Vol}(U) \choose j} j! ~
 \frac{(\mbox{Vol}(U)-j)! }{\left(\frac{\mbox{Vol}(U)-j}{2}\right)! ~ 2^{\frac{\mbox{Vol}(U)-j}{2}} }~
\frac{ (2E-\mbox{Vol}(U)-j)!}{\left(E-\frac{\mbox{Vol}(U)-j}{2}\right)! ~ 2^{\frac{2E-\mbox{Vol}(U)-j}{2}} }
\prod_{i=1}^{\nv} d_i!\:.
\end{equation}
}
The dependence of ${\cal H}$ on ${\bf d}$
has been suppressed for notational convenience. The first two terms on the 
right above count the number of ways we can choose 
$j$ configurations points each from $U$ and $U^c$ to match up. The
term $j!$ counts the number of ways of matching them. The remaining
configuration points have to be matched within the sets $U$ and $U^c$
as there are only $j$ edges crossing the cut. The number of ways of
doing this is the number of configurations on $U$ with $\mbox{Vol}(U)-j$ points,
times the number of configurations on $U^c$ with $2E-\mbox{Vol}(U)-j$ points, and 
with a degree sequence strictly bounded by ${\bf d}$ (since $j$ points 
each in $U$ and $U^c$ have been used up).
This yields the remaining terms in the bound above.
We obtain from (\ref{nconfig}) and (\ref{nconfigcut}) after some
simplification that
\begin{equation*}
\P(e_{U,U^c}(H) = j) = \frac{ {\cal H}_{U,U^c}(j) }{ H_{\nv}({\bf d}) } 
\le \frac{ {E\choose \mbox{Vol}(U)/2} {\mbox{Vol}(U)/2 \choose j/2 } {E-(\mbox{Vol}(U)/2) \choose j/2} }
{ {2E\choose \mbox{Vol}(U)} {j\choose j/2} } 2^j.
 \end{equation*}
Taking logarithms and using Stirling's formula, we get
\begin{eqnarray} 
\log \P(e_{U,U^c}(H) = j) &\le& 
E h \left( \frac{\mbox{Vol}(U)}{2E} \right) + \frac{\mbox{Vol}(U)}{2} h \left( \frac{j}{\mbox{Vol}(U)} \right)
+ \frac{2E-\mbox{Vol}(U)}{2} h \left( \frac{j}{2E-\mbox{Vol}(U)} \right) \nonumber \\
&& - 2E h \left( \frac{\mbox{Vol}(U)}{2E} \right) + O(\log \nv),
\label{stirling1}
\end{eqnarray}
where, for $x\in [0,1]$, $h(x) = -x\log x-(1-x)\log(1-x)$ is the binary
entropy of $x$.
Now, $2E = c\nv \log \nv$ and, since it was assumed that ${\bf d} \in
A_1(\alpha_1,\alpha_2)$, $|\mbox{Vol}(U) - cu\log \nv| \le u\sqrt{\alpha \log \nv}$, $\alpha=\max\{\alpha_1,\alpha_2\}$.
Moreover, $\tau \nv < u \le {\nv}/2$, while $j<\frac{1}{2}(1-\delta)c\nv\log\nv$. Hence, for some $\hat\delta_1$
and for large enough $\nv$, we have, for all $\delta\geq \hat\delta_1$
\[
h \left( \frac{j}{\mbox{Vol}(U)} \right) < 
h \left( \frac{(1-\delta)\nv\log\nv}{2\tau \nv \log \nv} \right)
= h\left(\frac{(1-\delta)}{2\tau}\right),
\]
and it can likewise be shown that, for some $\hat\delta_2$ and
for large enough $\nv$, we have, for all $\delta\geq \hat\delta_2$
$$h\left(\frac{j}{2E-\mbox{Vol}(U)}\right)< h(1-\delta)\:.$$ 
On the other hand, as $|U|<n/2$, for $n$ large,
\[
h \left( \frac{\mbox{Vol}(U)}{2E} \right) \geq h(\tau).
\]
Using the fact that $\mbox{Vol}(U)\leq 2E$, for all $U$, it follows from (\ref{stirling1}) that, for $\nv$ sufficiently large, 
\begin{equation} \label{stirling2}
\log \P(e_{U,U^c}(H) = j) \le -E \left( h\left(\tau\right)- h\left(\frac{1-\delta}{2\tau}\right)-h\left(1-\delta\right)\right)
\le -\kappa \nv \log \nv 
\end{equation}
where $\delta$ is chosen big enough so that $h(\tau)- h\left(\frac{1-\delta}{2\tau}\right)-h(1-\delta)>0$, i.e., $\kappa>0$.

The above bound applies for all subsets $U$ of $V$, of size $u$ where $\nv < u < \nv/2$. The number of subsets $U$ with cardinality between $\tau \nv$
and $\nv/2$ is smaller than the total number of subsets, which is $2^{\nv}$.
Hence, by the union bound,
\[
\P(H: \mbox{$\exists \, U$ with $\tau \nv < u < \nv/2$ and
$e_{U,U^c}(H)=j$}) \le 2^{\nv} e^{ -\kappa \nv \log \nv}.
\]
The above holds for each $j<\frac{1}{2}(1-\delta)cn\log n$. Applying the union bound once more,
\[
\P(H \in \hat {\cal E}_2(\nv,\delta,\tau,{\bf d})|{\bf d})
\le (1-\delta)cn\log{(n)}~ 2^{\nv-1} e^{ -\kappa\nv \log \nv }\:,
\]
for all ${\bf d} \in A_1(\alpha_1,\alpha_2)$.
Substituting this in (\ref{configbound2}) and noting that 
$\lambda=O(\log \nv)$, we see that, for $\delta$ large enough 
\[ 
\mu_{\nv}({\cal E}_2(\nv,\delta,\tau,{\bf d})|{\bf d}\in A_1(\alpha_1,\alpha_2)) \to 0 
\mbox{ as $\nv \to \infty$}\:.
\]

We also know from Theorem~\ref{thm:main} that
$\mu_{\nv}({\bf d} \notin A_1(\alpha_1,\alpha_2))$ goes to zero. Since
\[ \mu_{\nv}\left({\cal E}_2(\nv,\delta,\tau)\right) \le 
\mu_{\nv}\left({\cal E}_2(\nv,\delta,\tau,{\bf d})|{\bf d} \in A_1(\alpha_1,\alpha_2)\right)
+ \mu_{\nv}\left({\bf d} \notin A_1(\alpha_1,\alpha_2)\right),
\]
then, there exists $\delta _3>0$ such that $\mu_{\nv}({\cal E}_2(\nv,\delta_3,\tau)) \to 0$ as $\nv \to \infty$, 
as claimed.
\hfill $\Box$

Fix $\epsilon>0$ and $\tau <1/(1+4e)$, then by Propositions \ref{expander0}, \ref{expander1} and \ref{expander2}, there exists  $\tilde{\delta}$, independent of $n$,  which is the maximum of $\delta_1$, $\delta_2$ and $\delta_3$ for which the three propositions hold. Hence we have the following lower bound for the graph cut,
\begin{theorem} \label{expandermain}
For graphs $G$ drawn according to (\ref{eq:graph-distr}), there exists $\tilde{\delta}\in(0,1)$ such that for $U$ subset of $V$ with $u=|U|\leq n/2$, the number of edges crossing the cut $(U,U^c)$ is such that
\[
e_{U,U^c}\geq (1-\tilde\delta)  c u \log n, \quad \whp.
\]
\end{theorem}

\subsection{Conductance and Expansion}\label{sse-conductance}
Using Theorem \ref{expandermain}, we can easily recover asymptotic results on the conductance and the expansion of a graph drawn according to (\ref{eq:graph-distr}),  which are relevant for phenomena such as routing congestion analysis \cite{GMS03}, the behaviour of random walks in terms of the mixing and cover times \cite{lovasz96}, and epidemic threshold \cite{GMT05}.

Let $A=(a_{ij})_{i,j=1,\dots,n}$ be the adjacency matrix of a graph $G$ and $D=\mbox{Diag}(d_1,\dots,d_n)$ the diagonal matrix of the degree distribution of $G$.
First, we define the {\em isoperimetric constant or expansion} of a graph $G$  by
$$\phi=\inf_{U\subset V,\: u\leq
n/2 }\frac{e_{U,U^c}}{u}$$
It is related to $\lambda_2(L)$ the second (smallest) eigenvalue of the Laplacian $L=D-A$ of the graph through the following  inequality \cite{chung96,mohar}
$$\frac{\phi^2}{2d_{\max}}\leq \lambda_2(L)\leq 2 \phi\:.$$
The lower bound in the above inequality is known as the Cheeger's inequality.

The {\em conductance} of a graph $G$ is
defined by
$$\Phi=\inf_{U\subset V,\: \mbox{Vol}(U)\leq
 E }\frac{e_{U,U^c}}{\mbox{Vol}(U)}\:.$$

Let $\lambda_2(P)$ be the second (largest) eigenvalue of $P$ the transition matrix of the simple random walk on a graph $p_{ij}=a_{ij}/d_i$. By Cheeger's inequality \cite[Theorem 5.3]{lovasz96},

$$\frac{\Phi^2}{8}  \leq 1-\lambda_2(P)\leq \Phi\:.$$

\begin{theorem}\label{expansion conductance} For graphs $G$ drawn according to (\ref{eq:graph-distr}), and for the constant $\tilde{\delta}$ of Theorem \ref{expandermain}, we have that the expansion $\phi$ and the conductance $\Phi$ satisfy, 
$$(1-\tilde\delta) c\log{n} \leq \phi \leq c\log{n},\qquad (1-\tilde\delta) \leq \Phi \leq 1, \qquad \whp\:.$$
\end{theorem}
\noindent\emph{Proof}:
First note that if $d_{\min}$ is the minimum degree of $G$, then by Theorem \ref{thm:main}, $d_{\min}=c\log{n}-\sqrt{\alpha_1 \log{n}}$, $\whp$. Hence,
$$\phi\leq (1+o(1))c\log{n},\qquad \Phi \leq (1+o(1)), \qquad \whp\:.$$
The lower bounds follow from Theorem \ref{expandermain}.
\hfill $\Box$
\subsection{Failure resilience}\label{sse-resilience}

In the following, we work with graphs whose degree sequence belongs to the 
set $A_1(\alpha_1,\alpha_2)$ for some specified $\alpha_1$ and $\alpha_2$.
We are interested in the probability that the graph remains connected
when links fail independently with probability $p$. It is straightforward
to compute the probability that a given node $i$ becomes isolated due to
link failures; it is simply $p^{d_i}$. Thus, by the union bound, the
probability that some node becomes isolated is at most
\[
\sum_{i=1}^{\nv} p^{d_i} \le \nv p^{c\log \nv -\sqrt{\alpha_1 \log \nv}}
= \exp [ (1+c\log p)\log \nv -\sqrt{\alpha_1 \log \nv} \log p ].
\]
Hence, if $c\log p < -1$ or, equivalently, $p < \exp(-1/c)$, then the
probability that some node becomes isolated goes to zero as $\nv$
increases to infinity.

By way of comparison, consider the classical random graph model of 
Erd\"os and R\'enyi~\cite{ER60} with the same mean degree. Here, an
edge is present between each pair of nodes with probability
$c\log \nv/\nv$, independent of all other edges. Here we should assume that $c>1$ to ensure that the Erd\"os-R\'enyi graph is connected, $\whp$.
After taking failures
into account, the edge probability becomes $(1-p)c\log \nv/\nv$, and the
presence of edges continues to be mutually independent. It is well known
for this model that, if $(1-p)c< 1$, then the graph is disconnected with
high probability. Moreover, in a sense that can be made precise, the main 
reason for disconnection when $(1-p)c$ is ``close to" 1 is the isolation of
individual nodes. Intuitively, these arguments suggest that balanced
random graphs can tolerate link failure rates up to $e^{-1/c}$ while 
retaining connectivity, whereas classical random graphs can only tolerate 
failure rates up to $(c-1)/c$. We now rigourously establish a weaker result.

We shall use Thereom \ref{expandermain} to show that random graphs drawn
from the distribution $\mu_{\nv}$ can tolerate link failure rates up to
$\exp{\left(-\frac{1}{c(1-\tilde{\delta})}\right)}$, where $\tilde{\delta}$ is defined in Theorem \ref{expandermain}, without losing connectivity. 
\begin{theorem}
For
any $p<\exp{\left(-\frac{1}{c(1-\tilde{\delta})}\right)}$, a graph $G$ chosen at random from the distribution
$\mu_{\nv}$, and subjected to independent link failures with probability
$p$ remains connected, $\whp$.
\end{theorem}
\noindent\emph{Proof}:

Fix $p<\exp{\left(-\frac{1}{c(1-\tilde{\delta})}\right)}$. For a subset
$U$ of the vertex set, let $\hat e_{U,U^c}$ denote the number of
edges between $U$ and $U^c$ that have not failed. We shall show that,
with high probability, $\hat e_{U,U^c}>0$ for all subsets $U$, i.e.,
the graph is connected. Now, 
\[
\mu_n\left(\hat e_{U,U^c}(G)=0|e_{U,U^c}(G)\right) = p^{e_{U,U^c}(G)}.
\]

Assume that
$e_{U,U^c}(G) \ge (1-\tilde \delta) uc\log \nv$, for all $U\subseteq V$ 
with $u\leq \tau \nv$. Hence,
\[ 
\mu_n(\exists \, U: u\le \tau\nv, \, \hat e_{U,U^c}(G)=0)
\le \sum_{u=1}^{\tau\nv} {\nv \choose u} p^{(1-\tilde \delta) uc\log \nv}.
\]
Since $p<\exp{\left( -\frac{1}{c(1-\tilde{\delta})} \right) }$ given, then for some $\epsilon>0$ and $n$ large,  $p^{(1-\tilde\delta)c\log \nv}<e^{-(1+\epsilon) \log \nv}$. Using the inequality
${\nv \choose u} \le \nv^u/u!$, we get
\begin{eqnarray} \nonumber
\mu_n(\exists \, U: u\le \tau \nv, \, e_{U,U^c}(G)=0)
&\le & \sum_{u=1}^{\tau \nv} \frac{1}{u!} \left(np^{(1-\tilde\delta) c\log \nv} \right)^u\\ \label{disconnect-bd1}
&\le & \exp\left(np^{(1-\tilde\delta) c\log \nv}\right)-1 \leq \exp\left(ne^{-(1+\epsilon) \log \nv}\right)-1
\end{eqnarray}
which goes to zero as $\nv \to \infty$.

Suppose that $e_{U,U^c}(G) 
\ge (1-\tilde\delta) cu \log{n}$ for all $U\subseteq V$ with $\tau \nv <u \le \nv/2$. Hence,
\begin{equation} \label{disconnect-bd2}
\mu_n(\exists \, U: \tau \nv<u \le \nv/2, \, \hat e_{U,U^c}(G)=0)
\le \sum_{U:\tau \nv < u \le \nv/2} p^{(1-\tilde\delta)c u \log{n}}
\le 2^{\nv} p^{(1-\tilde\delta)\tau cn \log{n}}.
\end{equation}

We see from (\ref{disconnect-bd1}) and (\ref{disconnect-bd2}) that,
\[
\mu_n(\exists \, U: \hat e_{U,U^c}(G)=0 | e_{U,U^c}(G) 
\ge (1-\tilde\delta) cu \log{n}) \to 0
\mbox{ as } \nv \to \infty.
\]
Also, by Theorem \ref{expandermain},
\[
\mu_n( e_{U,U^c}(G) 
< (1-\tilde\delta) cu \log{n},\: \forall\: U\subseteq V, 0 <u \le \nv/2)
\to 0 \mbox{ as } \nv \to \infty,
\]
when $G$ is chosen according to the distribution $\mu_{\nv}$, which establishes the claim of the theorem.

\hfill $\Box$

\section{Appendix}\label{se-appendix}

Let $D_1,\ldots,D_{\nv}$ be iid random variables with distribution 
given by (\ref{eq:marginal}). Define
\begin{equation} \label{eq:defprob}
f(j,\gamma) = \frac{1}{j!} e^{-\beta j^2+\gamma j \log \nv}, \quad
\mbox{and} \quad F(\gamma) = \sum_{j=0}^{\infty} f(j,\gamma),
\end{equation}
so that $\PP(D_1=j) = f(j,\gamma)/F(\gamma)$. Now, the ratio
\[
\frac{f(j+1,\gamma)}{f(j,\gamma)} =
\frac{1}{j+1}e^{-(2j+1)\beta+\gamma \log \nv},
\]
is a decreasing function of $j$. Define $k_{\gamma}$ to be the smallest
value of $j$ for which $f(j+1,\gamma)/f(j,\gamma) \le 1$,
and note that the maximum of $f(j,\gamma)$ over $j$ is attained
at $k_{\gamma}$. Now, $k_{\gamma}-1$ is the integer part of the (unique)
solution of the equation
\begin{equation} \label{maxeq1}
h(x,\gamma) := - \log(x+1) - (2x+1)\beta + \gamma \log \nv = 0.
\end{equation}
It is readily verified that the solution is
\begin{equation} \label{argmax1}
x_{\gamma} = \frac{1}{2\beta} \Bigl( \gamma \log \nv + \log \log \nv
+ \frac{\gamma}{2\beta} \Bigr) + o(1).
\end{equation}
Let $k_{\gamma}=\left\lfloor x_{\gamma}\right\rfloor+1$. Then for any $j>0$,
\begin{eqnarray*}
\frac{f(k_{\gamma}+j+1,\gamma)}{f(k_{\gamma}+j,\gamma)} &=&
\frac{1}{k_{\gamma}+j+1} e^{-\beta(2k_{\gamma}+2j+1)
+\gamma \log \nv} \nonumber \\
&=& \frac{f(k_{\gamma}+1,\gamma)}{f(k_{\gamma},\gamma)}~
\frac{k_{\gamma}+1}{k_{\gamma}+j+1}~ e^{-2\beta j} \ \le \ e^{-2\beta j},
%\label{eq:ineqratio1}
\end{eqnarray*}
where we have used the fact that 
$f(k_{\gamma}+1,\gamma)/f(k_{\gamma},\gamma) \le 1$
to obtain the last inequality. Iterating this inequality yields
$ f(k_{\gamma}+j,\gamma)/f(k_{\gamma},\gamma) \le e^{-\beta j(j-1)}$. 
Similarly, we get
\begin{equation*} %\label{eq:ineqratio2}
\frac{f(k_{\gamma}-j-1,\gamma)}{f(k_{\gamma}-j,\gamma)} =
\frac{f(k_{\gamma}-1,\gamma)}{f(k_{\gamma},\gamma)} 
\Bigl( 1-\frac{j}{k_{\gamma}} \Bigr) e^{-2\beta j} \le e^{-2\beta j},
\end{equation*}
since $f(k_{\gamma},\gamma)/f(k_{\gamma}-1,\gamma) > 1$ by the
definition of $k_{\gamma}$. Iterating this inequality yields
$ f(k_{\gamma}-j,\gamma)/f(k_{\gamma},\gamma) \le e^{-\beta j(j-1)}$. Thus,
for all integers $j\ge -k_{\gamma}$, we have
the inequality
\begin{equation} \label{eq:ineqratio}
\frac{f(k_{\gamma}+j,\gamma)}{f(k_{\gamma},\gamma)} \le e^{-\beta |j|(|j|-1)}\le e^{-\beta(|j|-1)^2}.
\end{equation}

Next, we derive an equivalent for the above ratio.
Observe that, for any fixed $j$,
\begin{eqnarray*}
\frac{ f( k_{\gamma}+j, \gamma) }{ f( k_{\gamma}, \gamma) } &=&
\frac{k_{\gamma}!}{( k_{\gamma}+j )!}
e^{-\beta j(2k_{\gamma}+j)+\gamma j\log \nv} \\
&=& \frac{1}{k_{\gamma}^j} e^{-\beta j(2k_{\gamma}+j)+\gamma j\log \nv}
\Bigl( 1+O \Bigl( \frac{j^2}{k_{\gamma}} \Bigr) \Bigr).
\end{eqnarray*}
Taking logarithms, 
\begin{eqnarray*}
\log \frac{ f( k_{\gamma}+j, \gamma) }{ f( k_{\gamma}, \gamma) } &=&
-j\log k_{\gamma} -\beta j(2k_{\gamma}+j)+\gamma j\log \nv + 
O \Bigl( \frac{j^2}{\log \nv} \Bigr) \\
&=& j h(x_{\gamma},\gamma) + \alpha j - \beta j^2 +
O \Bigl( \frac{j^2}{\log \nv} \Bigr),
\end{eqnarray*}
where $\alpha = 2\beta (x_{\gamma}-k_{\gamma}+\frac{1}{2})$.
Note that $\alpha \in [-\beta,\beta]$ for all $\nv$ because $k_{\gamma}
\in [x_{\gamma},x_{\gamma}+1]$. Since $h(x_{\gamma},\gamma)=0$ by the
definition of $x_{\gamma}$, we can now write
\begin{equation} \label{eq:ratio}
g(j,\gamma) := \frac{ f( k_{\gamma}+j, \gamma) }{ f( k_{\gamma}, \gamma) }
= (1+\lambda_j) e^{\alpha j-\beta j^2}, \quad \mbox{where } \;
\lambda_j = O\Bigl( \frac{j^2}{\log \nv} \Bigr).
\end{equation}
Thus, by (\ref{eq:defprob}),
\begin{equation} \label{eq:Fgamma}
F(\gamma) = f(k_{\gamma},\gamma) \sum_{j=-k_{\gamma}}^{\infty}
g(j,\gamma) = K_0(\alpha,\beta) f(k_{\gamma},\gamma),
\end{equation}
where $K_0(\alpha,\beta) \sim \sum_{j=-\infty}^{\infty}
e^{\alpha j-\beta j^2}$ is bounded uniformly in $\gamma$ and $\nv$.

\subsection{Proof of Lemma \ref{le-deg dist1}}\label{proof-deg dist1}

We obtain from (\ref{eq:marginal}) and (\ref{eq:ratio}) that
\begin{eqnarray}
\EE D_1 &=&
\frac{ \sum_{j=0}^{\infty} jf(j,\gamma) }{ \sum_{j=0}^{\infty} f(j,\gamma) }
\ =\  k_{\gamma} 
\frac{ \sum_{j=-k_{\gamma}}^{\infty} ( 1+\frac{j}{k_{\gamma}} ) g(j,\gamma) }
{ \sum_{j=-k_{\gamma}}^{\infty} g(j,\gamma) } \nonumber \\
\noalign{\vspace{2mm}}
&=& k_{\gamma} \Bigl[ 1 + \frac{1}{k_{\gamma}}
\frac{ \sum_{j=-k_{\gamma}}^{\infty} j(1+\lambda_j) e^{\alpha j-\beta j^2} }
{ \sum_{j=-k_{\gamma}}^{\infty} (1+\lambda_j) e^{\alpha j-\beta j^2} } \Bigr] 
\nonumber \\
\noalign{\vspace{2mm}}
&=& k_{\gamma} + K_1(\alpha,\beta), \label{eq:mean}
\end{eqnarray}
where
\[
K_1(\alpha,\beta) \sim \Bigl( \sum_{j=-\infty}^{\infty}
je^{\alpha j-\beta j^2} \Bigr) \Bigm/ \Bigl( \sum_{j=-\infty}^{\infty}
e^{\alpha j-\beta j^2} \Bigr).
\]
Note that $K_1(\alpha,\beta)$ is bounded uniformly in $\gamma$ and $\nv$.
It is also easy to see that $\EE D_1$ is a continuous and
increasing function of $\gamma$. This yields the first claim of the lemma.

A similar calculation yields
\begin{eqnarray*}
\EE[(D_1)^2] &=&
\frac{ \sum_{j=0}^{\infty} j^2f(j,\gamma) }{ \sum_{j=0}^{\infty} f(j,\gamma) }
\ =\  k_{\gamma}^2 
\frac{ \sum_{j=-k_{\gamma}}^{\infty} (1+\frac{j}{k_{\gamma}})^2 g(j,\gamma) }
{ \sum_{j=-k_{\gamma}}^{\infty} g(j,\gamma) } \\
\noalign{\vspace{2mm}}
&=& k_{\gamma}^2 + 2k_{\gamma}K_1(\alpha,\beta)+K_2(\alpha,\beta),
\end{eqnarray*}
where
\[
K_2(\alpha,\beta) \sim \Bigl( \sum_{j=-\infty}^{\infty}
j^2 e^{\alpha j-\beta j^2} \Bigr) \Bigm/ \Bigl( \sum_{j=-\infty}^{\infty}
e^{\alpha j-\beta j^2} \Bigr)
\]
remains bounded, uniformly in $\gamma$ and $\nv$. Hence,
\[
\mbox{Var}(D_1) = K_2(\alpha,\beta) - K_1(\alpha,\beta)^2
\]
remains bounded. In fact, we see that $\mbox{Var}(D_1)$ is asymptotic
to the variance of a discrete Gaussian distribution; this distribution
is non-degenerate for any finite $\beta$. Hence, $\mbox{Var}(D_1)$
remains bounded below by some strictly positive constant as $\nv$ goes
to infinity. 

Next, we evaluate the moment generating function of $D_1$. 
Proceeding as in the calculations of the mean and variance, we have
\begin{eqnarray}
\EE[e^{\theta D_1}] &=&
\frac{ \sum_{j=0}^{\infty} e^{\theta j} f(j,\gamma) }{ \sum_{j=0}^{\infty}
f(j,\gamma) } \ =\  e^{\theta k_{\gamma}} \frac{
\sum_{j=-k_{\gamma}}^{\infty} e^{\theta j} g(j,\gamma) }{
\sum_{j=-k_{\gamma}}^{\infty} g(j,\gamma) } \nonumber \\
\noalign{\vspace{2mm}}
&=& e^{\theta k_{\gamma}} \frac{ \sum_{j=-k_{\gamma}}^{\infty}
(1+\lambda_j) e^{(\theta+\alpha) j-\beta j^2} }{ \sum_{j=-k_{\gamma}}^{\infty} 
(1+\lambda_j)e^{\alpha j-\beta j^2} } \nonumber \\
\noalign{\vspace{2mm}}
&\sim& e^{\theta k_{\gamma}} \frac{\psi(\theta+\alpha)}{\psi(\alpha)},
\label{eq:mgf1}
\end{eqnarray}
where 
\begin{equation} \label{eq:mgfgauss}
\psi(\theta) = \frac{\sum_{j=-\infty}^{\infty} e^{\theta j-\beta j^2}}
{\sum_{j=-\infty}^{\infty} e^{-\beta j^2}}
\end{equation}
is the moment generating function of the discrete Gaussian distribution 
which puts mass proportional to $e^{-\beta j^2}$ at each $j\in \Ints$.

\subsection{Proof of Lemma \ref{lem:tail1}}\label{proof-le-tail1}

We obtain using (\ref{eq:ineqratio}) and
(\ref{eq:Fgamma}) that, for $n$ large,

\begin{eqnarray*}
\PP(D_1> \nv^{1/4}) &=& \frac{ \sum_{j=n^{1/4}+1}^{\infty} f(j,\gamma) }
{ F(\gamma) } \\
&\le& \frac{1}{K_0(\alpha,\beta)} 
\sum_{j=0}^{\infty} e^{-\beta(j+ n^{1/4}-k_{\gamma})^2}\\
&\le& \frac{1}{K_0(\alpha,\beta)} 
\sum_{j=0}^{\infty} e^{-\beta(j+ \frac{1}{2}n^{1/4})^2}\\
&\leq &\frac{\sum_{j=0}^{\infty} e^{-\beta j^2}}{K_0(\alpha,\beta)}~ e^{-\beta\sqrt{n}/4}
\end{eqnarray*}
By the union bound 
\begin{equation*} 
\PP({\bf D} \in A_2^c) \le \sum_{i=1}^{\nv} 
\PP(D_i > {\nv}^{1/4}) \le K{\nv}e^{-\beta \sqrt{\nv}/4}.
\end{equation*}

which establishes the claim of the lemma. 

\subsection{Proof of Lemma \ref{lem:cgf}} \label{proof-le-cgf}

In what follows we prove the result for the sequence $D_i$. Following the same lines, one can prove the lemma for $\tilde D_i$.

Since $\EE D_1 = k_{\gamma}+K_1(\alpha,\beta)$, it follows from (\ref{eq:mgf1})
that
\begin{equation} \label{eq:mgf2}
\EE[e^{\theta X_{\nv 1}}] = e^{-\theta \EE[D_1]} \EE[e^{\theta D_1}]
\sim e^{-\theta K_1(\alpha,\beta)} \frac{\psi(\theta+\alpha)}{\psi(\alpha)}.
\end{equation}
For fixed $\theta$, this is bounded uniformly in $\nv$ since 
$K_1(\alpha,\beta)$ is so bounded, and $\psi$ does not depend on $\nv$. 
The first claim of the lemma now follows from the inequality 
$ \EE[e^{\theta |X_{\nv 1}|}] \le \EE[e^{\theta X_{\nv 1}}] + 
\EE[e^{-\theta X_{\nv 1}}]$.

Since $X_{\nv 1}=D_1-\EE D_1$, therefore $\mbox{Var}(X_{\nv 1})=\mbox{Var}(D_1)$,
and the second claim of the lemma is immediate from Lemma~\ref{le-deg dist1}.

The last claim of the lemma follows from the fact that
\begin{eqnarray*}
 \sum_{j=-\infty}^{\infty} \min \{ \PP(X_{\nv 1}=j), \PP(X_{\nv 1}=j+1) \} 
&= &\ \sum_{j=0}^{\infty} \min \{ \PP(D_1=j), \PP(D_1 = j+1) \} \\
&\ge &\ \sum_{j=0}^{\infty} \PP(D_1=j)\PP(D_1 = j+1) \\
&\sim & \ \frac{ \sum_{j=-\infty}^{\infty} (e^{\alpha j-\beta j^2})
(e^{\alpha(j+1)-\beta(j+1)^2}) }{ \sum_{j=-\infty}^{\infty}
e^{\alpha j-\beta j^2} } \ >\ 0.
\end{eqnarray*}

This completes the proof of the lemma. 

\subsection{Proof of Lemma \ref{lem:a1prob}} \label{proof:a1prob}

We shall bound $\PP({\bf D} \in A_1(\alpha_1,\alpha_2)^c)$ using
the moment generating function of $X_{\nv 1} := D_1-\EE D_1$, and 
Chernoff's bound. Observe from (\ref{eq:mgf2}) that
\begin{equation} \label{eq:mgf3}
\EE[e^{( \sqrt{\theta\log \nv})X_{\nv 1}}] = e^{- \sqrt{\theta\log \nv}
K_1(\alpha,\beta)} \frac{\psi( \sqrt{\theta\log \nv}+\alpha)}{\psi(\alpha)},
\end{equation}
where $\psi$ is defined in (\ref{eq:mgfgauss}). Here, $\alpha$ and $\beta$
are constants, and $K_1(\alpha,\beta)$ remains bounded as $\nv \to \infty$.
Let 
\[
y^* = \frac{ \sqrt{\theta\log \nv}+\alpha}{2\beta}, \quad
j^* = \lfloor y^* \rfloor.
\]
We have
\begin{eqnarray*}
&& \Bigl( \sum_{j=-\infty}^{\infty} e^{-\beta j^2} \Bigr) 
\psi( \sqrt{\theta \log \nv} +\alpha) \\
&&= \; e^{( \sqrt{\theta\log \nv}+\alpha)j^*-\beta (j^*)^2}
\sum_{j=-\infty}^{\infty} e^{( \sqrt{\theta\log \nv}+\alpha)(j-j^*)
-\beta(j^2-(j^*)^2)} \\
&&= \; e^{\beta j^* (2y^* - j^*)}
\sum_{k=-\infty}^{\infty} e^{2\beta (y^*-j^*)k -\beta k^2} \\
&&= \; e^{\beta (y^*)^2} e^{-\beta(y^*-j^*)^2}
\sum_{k=-\infty}^{\infty} e^{2\beta (y^*-j^*)k -\beta k^2}, 
\end{eqnarray*} 
and so,
\[
\psi( \sqrt{\theta\log \nv} +\alpha) = \kappa(\alpha,\beta,\theta) 
\exp \Bigl( \frac{( \sqrt{\theta\log \nv}+\alpha)^2}{4\beta} \Bigr),
\]
where $\kappa(\alpha,\beta,\theta)$ is bounded, uniformly in $\nv$ 
and $\theta$. Substituting this in (\ref{eq:mgf3}) yields
\begin{equation} \label{eq:mgf4}
\EE[e^{( \sqrt{\theta\log \nv})X_{\nv 1}}] = \kappa_1 \exp \Bigl(
\frac{\theta \log \nv}{4\beta} + \kappa_2  \sqrt{\theta\log \nv} \Bigr),
\end{equation}
where $\kappa_1$ and $\kappa_2$ may depend on $\alpha$, $\beta$, $\theta$
and $\nv$, but are bounded. Thus, we obtain using Chernoff's bound that
\[
\PP(X_{\nv 1} >  \sqrt{\alpha_2\log \nv}) \le \kappa_1 \exp \Bigl( 
-\sqrt{\theta \alpha_2} \log \nv + \frac{\theta \log \nv}{4\beta} + 
\kappa_2  \sqrt{\theta\log \nv} \Bigr),
\]
for all $\theta>0$. Take $\theta=4\alpha_2 \beta^2$. Now, by the union bound,
\[
\PP \Bigl( \bigcup_{j=1}^{\nv} \{ X_{\nv j} >  \sqrt{\alpha_2\log \nv} \} 
\Bigr) \le \kappa_1 \exp \left( -\left(\alpha_2 \beta-1\right) \log \nv + 
2\kappa_2  \beta \sqrt{\alpha_2\log \nv} \right).
\]
The constant $\alpha_2$ can be chosen large enough so that $\alpha_2 \beta-1 > K$. Hence the 
right hand side above decreases to zero faster than $e^{-K\log \nv}$
as $\nv \to \infty$. A similar bound can be obtained on the probability 
that $X_{\nv j} < - \sqrt{\alpha_1\log \nv}$ for some
$j\in \{ 1,\ldots,\nv \}$. Thus, we have shown that, given $K>0$,
we can choose $\tilde \alpha_1$ and $\tilde \alpha_2$ so that
\begin{equation} \label{eq:a1prob-bd1}
\PP \Bigl( \bigcup_{j=1}^{\nv} \{ X_{\nv j} >   \sqrt{\tilde\alpha_2\log \nv} 
\} \cup \bigcup_{j=1}^{\nv} \{ X_{\nv j} < -  \sqrt{\tilde\alpha_1\log \nv} 
\} \Bigr) < \frac{e^{-K\log \nv}}{2}
\end{equation}
for all $\nv$ sufficiently large. Here, $X_{\nv j} = D_j-\EE D_j$, and the
$D_j$ are iid with mean $c\log \nv$. Let ${\overline D}$ denote the
empirical mean of $D_1,\ldots,D_{\nv}$. The event, $|{\overline D}-\EE D_1|
>\sqrt{ \eta \log \nv}$ is the same as the event
$|X_{\nv 1}+\cdots+X_{\nv \nv}| > \nv \sqrt{\eta \log \nv}$. Using
the same Chernoff bound techniques as above, we can show that $\eta$
can be chosen so that, for sufficiently large $\nv$, this event has 
probability at most $e^{-K\log \nv}/2$. Combining this with 
(\ref{eq:a1prob-bd1}) yields the claim of the lemma: simply take 
$\sqrt{\alpha_1} = \sqrt{\tilde \alpha_1} + \sqrt{\eta}$ and $\sqrt{\alpha_2} = \sqrt{\tilde \alpha_2} + \sqrt{\eta}$.

\end{document}